# Virtual Target Selection for a Multiple-Pursuer Multiple-Evader Scenario

Isaac E. Weintraub[1] Alexander Von Moll[2] David Casbeer[3]
*Air Force Research Laboratory, Wright-Patterson AFB, OH, 45433, USA*

Satyanarayana G. Manyam[4]
*Infoscitex Corporation, Dayton, OH, 45431, USA*

**This paper considers an M-pursuer N-evader scenario involving virtual targets. The virtual targets serve as an intermediary target for the pursuers, allowing the pursuers to delay their final assignment to the evaders. However, upon reaching the virtual target, the pursuers must decide which evader to capture. It is assumed that there are more pursuers than evaders and that the pursuers are faster than the evaders. The objective is two-part: first, assign each pursuer to a virtual target and evader such that the pursuer team's energy is minimized, and second, choose the virtual targets' locations for this minimization problem. The approach taken is to consider the Apollonius geometry between each pursuer's virtual target location and each evader. Using the constructed Apollonius circles, the pursuer's travel distance and maneuver at a virtual target are obtained. These metrics serve as a gauge for the total energy required to capture a particular evader and are used to solve the joint virtual target selection and pursuer-evader assignment problem. This paper provides a mathematical definition of this problem, the solution approach taken, and an example.**

## I. Nomenclature

| | | |
|---|---|---|
| $c_{i,j}$ | = | The cost matrix for a pursuer $i$ and an evader $j$ |
| $E$ | = | Evader |
| $\mathcal{E}$ | = | Set of evaders |
| $\mathcal{H}$ | = | Hamiltonian |
| $\mathcal{H}_I$ | = | Hamiltonian before pursuer reaches virtual target |
| $\mathcal{H}_{\mathbb{I}}$ | = | Hamiltonian after the pursuer has reached the virtual target |
| $I$ | = | Interception point for $VT_i$ to evader $j$ |
| $J$ | = | Objective cost functional / performance measure |
| $L$ | = | Number of virtual targets |
| $M$ | = | Number of evaders |
| $M_V$ | = | Maximum number of virtual targets |
| $N$ | = | Number of pursuers |
| $O$ | = | Origin of Apollonius circle |
| $p$ | = | Costate |

[1]Electronics Engineer, Control Science Center, Air Force Research Laboratory, AIAA Senior Member.
[2]Aerospace Engineer, Air Force Research Laboratory, Control Science Center.
[3]Senior Electronics Engineer, Air Force Research Laboratory, Control Science Center, AIAA Associate Fellow.
[4]Research Scientist, Infoscitex Corporation.

| | | |
|---|---|---|
| $P$ | = | Pursuer |
| $\mathcal{P}$ | = | Set of pursuers |
| $R$ | = | Radius of the Apollonius circle |
| $t_0$ | = | Initial time |
| $t_1$ | = | Time when pursuer reaches virtual target |
| $t_f$ | = | Time time pursuer reaches evader at interception point |
| $v_E$ | = | Evader speed |
| $v_P$ | = | Pursuer speed |
| $\mathcal{V}$ | = | Set of virtual targets |
| $VT_i$ | = | Virtual Target $i$ |
| $\overline{VT}$ | = | Maximum bound on virtual target location |
| $\underline{VT}$ | = | Minimum bound on virtual target location |
| $x_{ij}$ | = | The selection matrix for a pursuer $i$ and an evader $j$ |
| $\lambda$ | = | Line of sight angle from virtual target to Apollonius circle origin |
| $\mu$ | = | Speed ratio of the evader with respect to the pursuer |
| $\sigma_E$ | = | Evader angle from Apollonius circle origin to interception location |
| $\sigma_P$ | = | Pursuer angle from Apollonius circle origin to interception location |
| $\psi_E$ | = | Heading of the evader in the global fixed frame |
| $\psi_P$ | = | Heading of the pursuer in the global fixed frame |

## II. Introduction

Pursuit assignment problems involving mobile agents represent a relevant and popular class of problems for the aerospace and defense community, amongst others. In this problem, mobile vehicles or agents aim to "capture" other mobile or static targets through collocation or reaching the targets within a non-zero range. The task of assigning pursuers to targets is an area of research commonly referred to as "target assignment." The interested reader may want to read some works that review and survey the problem of weapon target assignment [1],[2].

In this paper a set of mobile agents, called pursuers, aim to capture a set of constant-speed fixed-course targets, called "evaders."" While the evaders in this work do not actively evade capture by the pursuers, the term evader is used to align with the context of differential games [3]. This paper focuses solely upon the collective strategies of the pursuer team assuming a known strategy for the evader team. Furthermore, the term "virtual target" is used to describe an intermediate location whereby the pursuers make a decision as to which evader to engage in pursuit. Using the term "evader" distinguishes between the virtual targets (that do not move) and evaders (that are fixed in speed and course).

What delineates this work from many available target assignment problems in literature is that the pursuer team first navigates to a set of virtual target locations whereby the pursuers decide which evader to pursue. This allows the pursuer team an opportunity to delay their guidance solutions in the event that some evaders may be captured prior to the delayed decision of the pursuers. As a result, the virtual targets allow the assignment of pursuers to evaders to occur later in time than at the onset. This work aims to address the question as to where the virtual targets should be located so as to minimize the overall pursuer team's energy to reach the evaders. The energy is modeled as the range traversed by each of the pursuers as well as the maneuver made by the pursuers upon reaching the virtual target.

### A. Multi-Pursuer Multi-Evader Games

Multi-pursuer-multi-evader problems involve a conflict between two teams. The team of pursuers aims to generally minimize some desired objective/performance measure while the team of evaders aims to maximize the desired objective/performance measure. In capture games, the pursuers aim to capture the evaders while the evaders strive to evade capture. In target defense and border defense games, the evaders aim to reach a target location or region while the pursuers strive to prevent the evaders from doing so – in this case the evaders are commonly called attackers and the pursuers are commonly called defenders. The outcome and optimal strategies for the pursuer and evader teams have been investigated in the context of differential games [3]; in the context of this paper, we consider the evader's strategies to be fixed and known to the pursuers. The single-sided optimization that provides the optimal pursuer strategies, provided the evaders take a fixed course, can lead to differential game solutions and is left for future work. The aim of this paper is to consider only the pursuer's strategies and investigate the role of the virtual targets in the pursuit scenario.

Many related works have considered multi-agent pursuit-evasion. In the references that follow, the relation to the work described herein is in how the assignment of pursuers to evaders is made, what objective/performance measures are optimized, and how pursuers cooperate. With this in mind, the following is a description of various related works in the field of multi-agent pursuit-evasion games. In this vein, the interested reader is encouraged to read a nice review of multi-agent pursuit-evasion perimeter defense games that can be found in a work by Shishika and Kumar [4]. In a very related work by Michael Day [5], the task assignment of multiple pursuers against multiple evaders is considered. Day considered pursuers to have a multi-phased strategy for capturing evaders including an ingress phase, assignment phase, and pursuit phase. Katz et al. describe a zero-sum differential game formulation whereby a red force tries to invade a blue territory and the blue force aims to prevent red's invasion [6]. Their formulation also incorporated avoidance of obstacles and used a quadratic cost formulation. Stipanović, Melikyan, and Hovakimyan considered non-holonomic kinematic model for a multi-pursuer multi-evader scenario in the plane, and leveraged a Lyapunov approach opposed to solving the classical Hamilton-Jacobi-Isaacs partial differential equations for a value function [7]. The methods are implemented in state-feedback form and they demonstrated that the approach taken is an approximation of the equilibrium strategies.

### B. Multi-Pursuer Multi-Evader Border/Target Defense Games

Multi-agent pursuit-evasion games involving border defense and target defense consider the defense of a region or specified mobile or static target/agent. Cruz et al. posed a two-team game composed of heterogeneous agents wherein an attacking team strives to reach some fixed targets and a defending team aims to defend the targets by capturing the attackers prematurely [8]. In this work non-zero sum dynamic games were conducted where agents implement the one-time Nash optimizing strategies. In a work by Earl and D'Andrea a team of defenders aim to protect a region of space from pursuing attackers [9]. This game was posed and solved using a mixed-integer linear program solver – the optimized strategies for defenders and attackers was performed at every time step and the authors comment that this is a computationally intensive task for the static environment. Rusnak posed a multiplayer differential game called, "The lady, the bandits, and the body-guards – a two team dynamic game" [10]. In this paper a two team dynamic game consists of bandits who aim to capture a lady while the lady, teamed with body-guards aims to prevent capture of the lady. Rusnak's work entails the pursuing team trained on a specific high-value target while a defensive team aims to protect the high-value target. In a very related work, Garcia et al. pose a multi-pursuer multi-evader differential game wherein $N$-pursuers and $M$-evaders are considered [11]. The evaders aim to reach a boundary while evading the pursuers and the pursuing team aims to capture the evaders before they are able to reach a boundary. The resulting strategies leverage Apollonius circle geometry to compute all possible pursuer-evader assignments. Using all possible assignments, the strategies for the evaders and pursuers are obtained. Yan et al. posed a 3-D heterogeneous reach avoid differential game involving heterogeneous players with different speeds (pursuers faster than evaders) and capture radii [12]. One team strives to prevent the other from reaching a plane and the pursuers are assigned to evaders using a bipartite graph. Garcia et al. considered the defense of a plane in 3-D with only one evader and a team of pursuers, demonstrating that three pursuers were required to ensure successful defense [13]. In their work point-capture was considered and both pursuers and evaders had the same speed. In a very relevant work by Asgharnia, Schwartz, and Atia, a reinforcement learning approach is taken to find optimal attacker and defender strategies for guarding a target where defenders are able to change their policy mid-game [14]. Zepp, Fu, and Liu posed a target defense scenario wherein a superior defender captures a team of attackers in sequence to protect a target [15], [16]. Fu and Liu also considered a team of slower defenders against a superior pursuer [17].

### C. Multi-Pursuer Multi-Evader Games of Incomplete Information

Games involving multiple pursuers and evaders with incomplete information represent a class of games where agents have limited or acquired information about their adversary. Antoniades, Kim, and Sastry considered a multi-pursuer-multi-evader game of incomplete information in [18]. In their work the players operate in a bounded region and perform heuristic strategies – the pursuers minimize sensor overlap between them and the evaders strive to evade capture by the pursuers. Li, Cruz, Chen, Kwan, and Chang considered a deterministic multi-player pursuit-evasion differential game in [19]. In their paper the pursuer-evader assignment problem leverages a hierarchical approach that considered a combinatorial approach to leveraging single-pursuer single-evader engagements in a deterministic multi-player game. Later, Li, Cruz, and Schumacher considered multi-player stochastic pursuit-evasion differential games involving multiple pursuers and evaders [20]. In their work pursuers were endowed with a non-zero capture range and a limited look-ahead was provided to agents. Wei et al. proposed a decentralized approach to multiplayer pursuit-evasion games where agents are endowed with jamming capabilities [21] – the estimation error of the opposing agents grows. The paper presents a unique game of information and capture/evasion under the pretense that agent actions are based upon estimates of their adversary.

### D. Multiplayer Reach-Avoid Games

Multi-player reach avoid games aim to obtain the capability of teams of agents by providing bounds on performance for pursuers and or evaders. In general, these games involve pursuit-evasion amongst obstacles and establish reachable regions for pursuers and or evaders. Chen, Zhou and Tomlin investigated a multi-player reach avoid target-defense game wherein $N_a$ attackers and $N_d$ defenders play on a compact domain with obstacles [22]. In their paper they provide guarantees for an upper bound, leveraging a bipartite graph structure, on the number of attackers that are able to reach a designated target. In another work, Zhou et al. Consider a similar scenario, but where one team's actions are provided at the onset [23].

### E. Multiple Pursuer Single Evader

In many cases of multi-pursuer multi-evader games, the roles of pursuers are extended from the cases where there exists only one evader. These single-evader problems can either involve a superior evader (more maneuverable or faster) or are used for extending optimal multi-pursuer single-evader strategies to a class of multi-pursuer multi-evader scenarios. While our paper addresses the general multi-pursuer multi-evader scenario; a degenerate case that only involves one evader is also possible. For this reason, related work in multi-pursuer single-evader is presented as it relates to this work.

Awheda and Schwartz considered a multi-pursuer single-superior-evader pursuit-evasion differential game [24]. In their formulation, they leveraged the Apollonius circle as a mechanism to define capture regions of each pursuer with respect to the evader. Using this information a fuzzy-based distributed approach was used to trap and capture the superior evader. Zhou et al. considered a pursuit-evasion game wherein a number of pursuers attempt to capture a single evader. Their work presents a decentralized real-time algorithm for cooperative pursuer and leverages Voronoi partitions [25]. Bakolas and Tsiotras have considered using Voronoi partitions for assigning pursuers to a static target [26], and to a mobile target [27]. Makkapati and Tsiotras leveraged Apollonius circle geometry in order to construct a Voronoi partition used to dictate the assignment of pursuers to evaders [28]. In their work the partitioning provided min-time capture strategies for the pursuers and partitions were recomputed after each individual evader was captured. Ibragimov et al. have also considered a unique perspective at pursuit of a single evader [29]. In their work they consider a differential game of countably many pursuers against a single evader and bound the terminal time – an investigation of the game of kind is performed whereby the optimal lower bound of the distance between the evader and the pursuers. A more recent work by Ibragimov considered a multi-pursuer single evader differential game under Gronwall constraints [30]. Garcia considered the containment of a high-speed evader by multiple slow pursuers [31]. Also included in that work were pursuers modeled to have non-zero capture radius. Garcia's work is related to a work by Jin and Qu that analyzed the containment of a superior evader by a number of pursuers [32]. In their work, pursuers leveraged point-capture and the methodologies leveraged Apollonius circle geometry for obtaining optimal strategies. Chen et al. considered the game of kind rather than the game of degree, investigating the initial conditions for which the evader was successful in avoiding capture or where pursuers' capture of the evader was guaranteed [33].

### F. Delayed Decision Guidance and Virtual Targets

This paper makes use of virtual targets to delay the pursuer-evader assignment. The delay of decision-making provides the ability for pursuers to be robust to changes during a general multi-pursuer multi-target engagement scenario [34]. In the work by Turetsky, Weiss, and Shima, it was assumed that the pursuer, as posed, was part of a salvo of multiple pursuers, and that a number of targets may have been captured before the pursuers have assigned which targets to engage. In their work, pursuers navigate toward the targets and delay their decisions, selecting highest priority targets. The end result is a more robust strategy for ensuring capture of incoming targets and preventing pursuers from being unassigned. Weiss, Shalumov and Shima later investigated how the delayed decision guidance strategies scale in a *many-on-many* intercept scenario [35]. In this work the pursuers navigate to static virtual target locations prior to assigning their targets. Merkolov, Weiss, and Shima leveraged the use of virtual targets for emulating more advanced guidance laws in [36] by making the virtual targets previously seen in [35] to be moving targets. These works leverage the use of virtual targets to improve the performance of target capture; but, the question remains as to the location of virtual targets that minimize the energy requirements of the pursuer team thereby maximizing the reachable range and g's available at interception – this paper addresses this open problem.

### G. Combinatorics

In weapon target assignment problems, a set of weapons must be assigned to a set of targets. For every pair of a weapon and a target, the intercept time or energy expenditure could be the cost associated with each assignment; and the minimum cost assignment can be posed as a bipartite matching problem [37]. In the current paper, we assume the targets travel on a fixed course at constant speed, and using elementary geometry, we can evaluate the interception position starting from a virtual target location. However, the virtual target location (within some predefined region) for every pursuer must also be determined. We sample the feasible region and refer to this discrete set of potential locations for virtual targets as $\mathcal{V}$. For a given pursuer-virtual target-evader combination, we can compute the cost of interception. The combinatorial problem is to find the virtual target and evader assignment for every pursuer, such that the total assignment cost is minimized. The resulting problem is a variant of the tripartite matching problem or three index assignment problem [38], [39], [40], with an additional constraint for the maximum number of virtual targets selected. We formulate this as a mixed integer linear program and solve using a branch and cut framework. The solution selects the virtual target locations from a discrete set of points sampled from a candidate in a specified region, and gives the assignments that minimizes the total cost.

### H. Paper Organization

This paper is presented as follows: In Section III, the problem definition and general pursuer-virtual-target-evader geometry is provided. Also, the performance measure for the pursuer is mathematically defined, describing how the energy of the pursuer can be obtained through simple geometry. In Section IV, the task assignment of pursuers to virtual targets and then to evaders is described as a linear program. In Section V, an example scenario is investigated and the results are presented. Finally, in Section VI, final concluding remarks and extensions are presented.

## III. Problem Definition

In order to identify the location of the virtual targets that minimize the energy of a pursuer to capture a non-maneuvering evader. While the strategic purpose of these virtual targets offers some known benefits [35], the question remains as to the location of these targets that can serve to minimize the time of arrival of the pursuer to the evader and also minimize the maneuver made when the pursuer reaches the virtual target. A generalized figure that highlights this engagement can be seen in Fig. 1.

Consider the $N$-pursuer, $M$-evader, $L$-virtual target scenario. That is there are many pursuers, evaders, and virtual targets. It is assumed that all of the pursuers are faster than all of the evaders: $v_{P_i} > v_{E,j} \forall i \in [1..N], j \in [1..M]$. This ensures that every pursuer is capable of capturing any of the evaders; relaxing this restriction is left as future work. Also, it is assumed that capture is attained when a pursuer collocates with an evader – zero capture-range is assumed. While, in a general sense there are many pursuers, evaders, and virtual targets, for the sake of providing non-trivial and invalid examples, it is assumed that the number of pursuers is greater than or equal to the number of evaders, $N \geq M$ – capture of all evaders is guaranteed. It is also assumed that the number of virtual targets is less than the number of pursuers. Provided these assumptions it can be stated that $N > L$. Furthermore, it is assumed that all agents have unlimited range of movement – the concept of range-limited pursuit-evasion is not considered in this work.

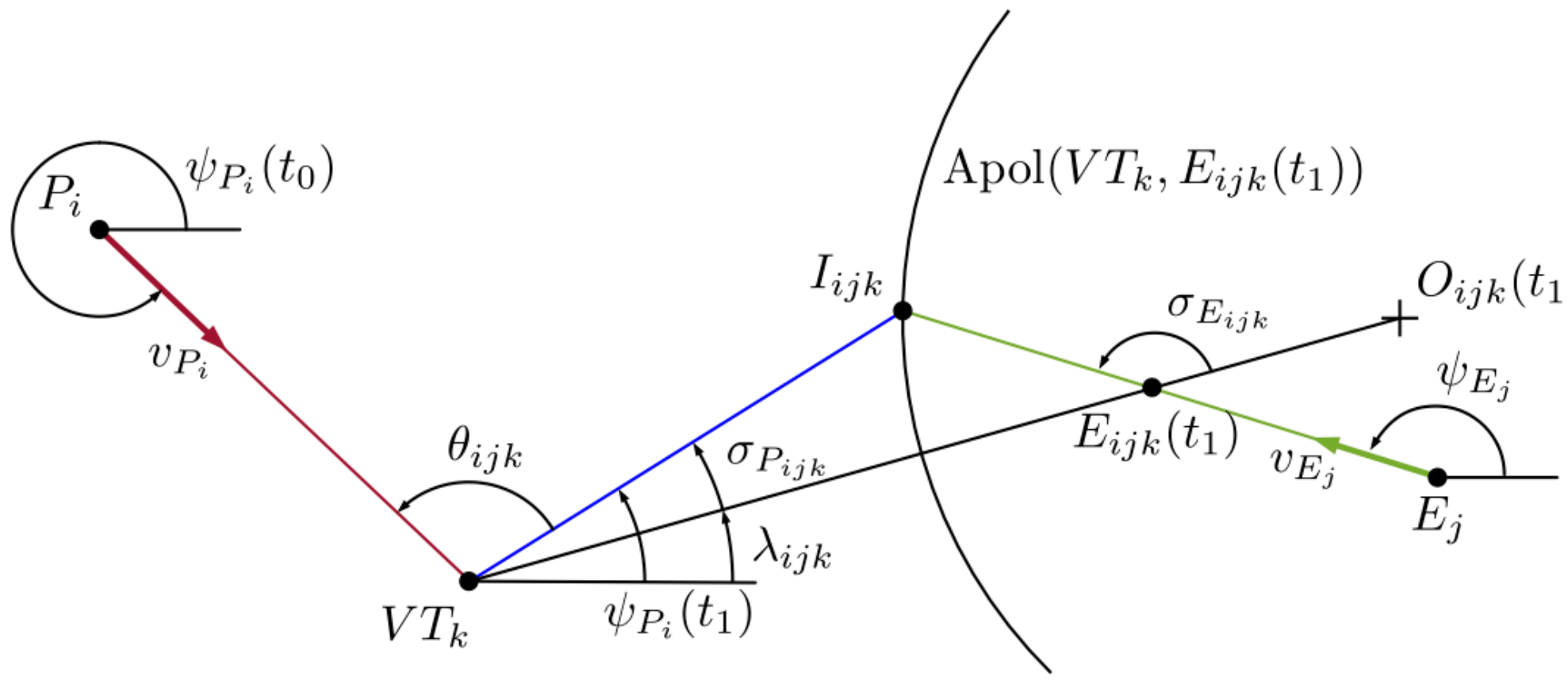

**Fig. 1 The generalized geometry of a pursuer, $P_i$, navigating to a virtual target, $VT_k$, and then choosing a heading that intercepts a non-maneuvering Evader, $E_j$, intercepting at the location $I_{ijk}$.**

Provided these assumptions, the location of the pursuers and evaders is assumed to be known at the onset. Since the evaders are non-maneuvering, they do not require knowledge of the pursuers or the location of virtual targets; also, their strategies are fixed and known to the pursuers. The pursuers are assumed to be able to share information so as to assign their selection of evaders and virtual targets accordingly. It is also assumed that the pursuers do not collide with one another, even if they select the same virtual targets; deconfliction of pursuers is left as future work. The task is to obtain optimal virtual target locations such that the overall energy of the pursuing team is minimized.

Consider the motion of $N$ pursuers and $M$ evaders as holonomic vehicles:

$$\begin{aligned}
\dot{x}_{P_i} &= v_{P_i} \cos \psi_{P_i} \forall i \in [1...N],\\
\dot{y}_{P_i} &= v_{P_i} \sin \psi_{P_i} \forall i \in [1...N],\\
\dot{x}_{E_j} &= v_{E_j} \cos \psi_{E_j} \forall j \in [1...M],\\
\dot{y}_{E_j} &= v_{E_j} \sin \psi_{E_j} \forall j \in [1...M].
\end{aligned}$$

The position of the $N$ pursuers are $P_i = (x_{P_i}, y_{P_i}) \in \mathbb{R}^2$, $P \in \mathbb{R}^{N \times 2}$; the position of the $M$ evaders are $E_j = (x_{E_j}, y_{E_j}) \in \mathbb{R}^2$, $E \in \mathbb{R}^{M \times 2}$; and the location of the $L$ virtual targets are $VT_k = (x_{VT_k}, y_{VT_k}) \in \mathbb{R}^2$, $VT \in \mathbb{R}^{L \times 2}$. The heading for the evaders $\psi_E \in \Theta^M, \Theta \in [0, 2\pi] \subset \mathbb{R}$ is constant and known by the pursuers. The control for the pursuers are the instantaneous headings for the pursuers, $\psi_P \in \Theta^N$.

Given a particular virtual target and evader, the objective of an individual pursuer is to reach the evader in minimum time and with minimum maneuver at the virtual target. The overall objective to be minimized is

$$\min_{\psi_P} J = \min_{\psi_P} \left\{ \pi - \theta_{ijk} + \int_{t_0}^{t_f} 1 dt \right\} = \min_{\psi_P} \left\{ \pi - \theta_{ijk} + \int_{t_0}^{t_1} 1 dt + \int_{t_1}^{t_f} 1 dt \right\}, \tag{1}$$

where the maneuver made by the pursuer at the virtual target is $\theta_{ijk}$, the time $t_1$ is when the pursuer $P_i$ reaches the assigned virtual target $VT_k$, and the final time $t_f$ is when the pursuer, $P_i$ is collocated with the assigned evader $E_j$. As will be seen later, the maneuver made by the pursuer when exiting the virtual target location will be represented by the angle $\pi - \theta_{ijk}$ and will be defined after we show that the time-optimal paths are, in fact, straight lines. The reason for the subscript $ijk$ is to annotate that this maneuver is specific for a specific pursuer, virtual target, and evader selection; if any of these designated agents or virtual targets change; the angle will invariably be different.

Each of the three components comprising the objective function (Eq. (1)) are described in the following subsections.

### A. Pursuer to Virtual Target Min-Time Strategy

The objective is to find the minimum time path of a pursuer $P_i$ to navigate to a desired virtual target $VT_k$. The objective cost for this phase is:

$$\min_{\psi_{P_i}} J = \min_{\psi_P} \left\{ \int_{t_0}^{t_1} 1dt \right\} = \min_{\psi_P} \{t_1 - t_0\}. \tag{2}$$

Assuming that the initial time $t_0 = 0$ in Eq. (2), the objective is to minimize the time $t_1$ wherein the pursuer arrives at the virtual target location. The Hamiltonian:

$$\mathcal{H}_I = p_{x_{P_i}} v_{P_i} \cos \psi_{P_i} + p_{y_{P_i}} v_{P_i} \sin \psi_{P_i}, \tag{3}$$

By optimal control theory, the necessary conditions for optimality provide the optimal costate dynamics as:

$$\dot{p}_{x_{P_i}} = -\frac{\partial \mathcal{H}_I}{\partial x_{P_i}} = 0, \quad \dot{p}_{y_{P_i}} = -\frac{\partial \mathcal{H}_I}{\partial y_{P_i}} = 0 \tag{4}$$

Thus, the optimal costates are constant to reach the virtual target in minimum time. Using the stationarity equation:

$$\frac{\partial \mathcal{H}_I}{\partial \psi_{P_i}} = 0 \Rightarrow -p_{x_{P_i}} v_{P_i} \sin \psi_{P_i} + p_{y_{P_i}} v_{P_i} \cos \psi_{P_i} = 0 \tag{5}$$

The optimal control is therefore constant since it solely relies upon the costates which are constant. Therefore the optimal min-time path for the pursuer to reach the virtual target is a straight-line path.

Considering the navigation of a pursuer to a virtual target, the time-optimal path from a pursuer, $P_i$ to a virtual target $VT_k$, is a straight line. Therefore the optimal heading that the pursuer takes is:

$$\psi^*_{P_i}(t) = \left\{ \text{atan2}\left( y_{VT_k} - y_{P_i}, x_{VT_k} - x_{P_i} \right) \mid t \in [t_0, t_1) \right\} \tag{6}$$

### B. Pursuer (From Virtual Target) to Evader Min-Time Strategy

The objective is to find the minimum time path of a pursuer $P_i$ (currently located at the virtual target) to navigate to a desired Evader $E_j$. This is the second phase of the pursuer's path. The objective is:

$$\min_{\psi_{P_i}} J = \min_{\psi_{P_i}} \left\{ \int_{t_1}^{t_f} 1dt \right\} = \min_{\psi_{P_i}} \{t_f - t_1\} \tag{7}$$

Similar to the first phase, the Hamiltonian for the second phase is:

$$\mathcal{H}_{\mathbb{I}} = p_{x_{P_i}} v_{P_i} \cos \psi_{P_i} + p_{y_{P_i}} v_{P_i} \sin \psi_{P_i} \tag{8}$$

following the same procedure as performed in the first phase, the necessary conditions for optimality provide that the costates are stationary,

$$\dot{p}_{x_{P_i}} = -\frac{\partial \mathcal{H}_{\mathbb{I}}}{\partial x_{P_i}} = 0, \quad \dot{p}_{y_{P_i}} = -\frac{\partial \mathcal{H}_{\mathbb{I}}}{\partial y_{P_i}} = 0 \tag{9}$$

Considering the stationarity condition for the second phase:

$$\frac{\partial \mathcal{H}_{\mathbb{I}}}{\partial \psi_{P_i}} = 0 \Rightarrow -p_{x_{P_i}} v_{P_i} \sin \psi_{P_i} + p_{y_{P_i}} v_{P_i} \cos \psi_{P_i} = 0. \tag{10}$$

The optimal control is therefore constant since it solely depends upon the costates which are stationary. Therefore the optimal min-time path for the pursuer to reach the evader is a straight-line path.

Because the path of the pursuer is a straight line and it is aiming to intercept a slower moving evader that is also moving in a straight line path, the use of Apollonius circle is a useful tool [41]. Apollonius circle provides the locus of possible interception points for two non-maneuvering agents of different speeds. In order to construct the Apollonius circle for the pursuer and evader, the location of both agents (foci), their speed ratio between the evader and the pursuer must be known. Since the heading of the evader is known, the interception point of the evader by the pursuer is obtained through the intersection of the evader's path with the Apollonius circle.

Since the evader was moving while the pursuer was navigating to the virtual target, the foci of the Apollonius circle requires forward propagating the evader's location to a place in the Cartesian space when the pursuer arrives at the virtual target. The time that occurs is $t_1$ and it is obtained as follows:

$$t_1 = \frac{\sqrt{\left(x_{VT_k} - x_{P_i}(t_0)\right)^2 + \left(y_{VT_k} - y_{P_i}(t_0)\right)^2}}{v_{P_i}} \tag{11}$$

The evader's location when the pursuer reaches the virtual target $E_{ijk}(t_1)$ is therefore:

$$E_{ijk}(t_1) = E_j(t_0) + t_1 v_{E_j} \begin{pmatrix} \cos \psi_{E_j} \\ \sin \psi_{E_j} \end{pmatrix} \tag{12}$$

The line of sight angle from the virtual target $VT_k$ to the evader $E_{ijk}(t_1)$ when the pursuer is at the virtual target is

$$\lambda_{ijk} = \operatorname{atan2}\left(y_{E_{ijk}}(t_1) - y_{VT_k}, x_{E_{ijk}}(t_1) - x_{VT_k}\right) \tag{13}$$

The distance between the virtual target and the evader when the pursuer is at the virtual target is:

$$\overline{VT_k E_{ijk}} = \sqrt{\left(x_{E_{ijk}}(t_1) - x_{VT_k}\right)^2 + \left(y_{E_{ijk}}(t_1) - y_{VT_k}\right)^2} \tag{14}$$

The speed ratio between the pursuer and the evader is:

$$\mu_{ij} = \frac{v_{E_j}}{v_{P_i}} \tag{15}$$

The origin of the Apollonius circle lies on the ray from the pursuer to the evader, now $VT_k$ and $E_{ijk}(t_1)$. The offset from the evader to the origin, opposite the direction of the pursuer is:

$$\overline{E_{ijk}(t_1) O_{ijk}(t_1)} = \frac{\mu_{ij} \overline{VT_k E_{ijk}(t_1)}^2}{1 - \mu_{ij}^2} \tag{16}$$

and therefore the origin of the Apollonius circle when the pursuer is at the virtual target is

$$O_{ijk}(t_1) = E_{ijk}(t_1) + \overline{E_{ijk}(t_1) O_{ijk}(t_1)} \begin{pmatrix} \cos \lambda_{ijk} \\ \sin \lambda_{ijk} \end{pmatrix} \tag{17}$$

The radius of the Apollonius circle is

$$R_{ijk} = \overline{O_{ijk}(t_1) E_{ijk}(t_f)} = \overline{O_{ijk}(t_1) I_{ijk}} = \frac{\mu_{ij} \overline{VT_k E_{ijk}(t_1)}}{1 - \mu_{ij}^2} \tag{18}$$

Using the heading of the evader in the global fixed frame, the local angle of the evader along the line of sight as shown in Fig. 1 is $\sigma_{E_{ijk}}$ and is computed by subtracting the line of sight angle from the evader's heading:

$$\sigma_{E_{ijk}} = \psi_{E_j} - \lambda_{ijk} \tag{19}$$

Using this local angle and the geometry of Apollonius circle, the resulting angle of the pursuer along the line of sight angle, $\sigma_{P_{ijk}}$ as shown inFig. 1 is

$$\sigma_{P_{ijk}} = \sin^{-1}\left(\mu_{ij} \sin \sigma_{E_{ijk}}\right) \tag{20}$$

To obtained the heading of the pursuer that reaches the evader in minimum time, we add the line of sight angle to find

$$\psi^*_{P_i}(t) = \left\{\sigma_{P_{ijk}} + \lambda_{ijk} \mid t \in \left[t_1, t_f\right]\right\} \tag{21}$$

From the geometry of Apollonius circle, distance traversed by the pursuer is

$$\overline{VT_k I_{ijk}} = \frac{\overline{VT_k E_{ijk}}}{1-\mu_{ij}^2}\left(\mu_{ij}\cos\sigma_{E_{ijk}} + \sqrt{1-\mu_{ij}^2\sin^2\sigma_{E_{ijk}}}\right) \tag{22}$$

The resulting interception point is therefore:

$$I_{ijk} = \overline{VT_k I_{ijk}}\begin{pmatrix}\cos\left(\sigma_{P_{ijk}} + \lambda_{ijk}\right)\\ \sin\left(\sigma_{P_{ijk}} + \lambda_{ijk}\right)\end{pmatrix} + \begin{pmatrix}x_{VT_k}\\ y_{VT_k}\end{pmatrix} \tag{23}$$

### C. Maneuver Required at Virtual Target

Because the min-time optimal path for the pursuer for both phase-I and phase-II are straight line paths, there exists an interior angle between these two locations. This interior angle represents the maneuver required by the pursuer to reach the evader in minimum time. It is the aim to not only minimize the time the pursuer takes to reach the evader, but also minimize the maneuver required at the virtual target. Using the law of cosines, the interior angle is

$$\theta_{ijk} = \cos^{-1}\left(\frac{\left\langle\overrightarrow{P_i VT_k}, \overrightarrow{VT_k I_{ijk}}\right\rangle}{\sqrt{\left\langle\overrightarrow{P_i VT_k}, \overrightarrow{P_i VT_k}\right\rangle}\sqrt{\left\langle\overrightarrow{VT_k I_{ijk}}, \overrightarrow{VT_k I_{ijk}}\right\rangle}}\right) \tag{24}$$

Since the objective is to minimize the maneuver required, the task is to make $\theta_{ijk}$ as close to $\pi$ as possible. Since the range of $\cos^{-1}(\cdot) \in [0, \pi]$, the cost of $\pi - \theta_{ijk}$ is selected. By minimizing the cost of $\pi - \theta_{ijk}$, the min-maneuvering path would be one wherein $\theta_{ijk} = \pi$.

## IV. Task Assignment

The task at hand is to minimize the time it takes for the pursuer to reach its assigned evader by way of an assigned virtual target as well as the maneuver required at the virtual target. The position of the virtual targets could be selected from a closed continuous set of locations, however the number of locations for virtual targets is limited to be less than prescribed maximum $M_V$. The assignment problem aims to select the pursuer-virtual target-evader combinations for every pursuer such that the sum of costs from all assignments is minimized. To pose a tractable problem, we sample the possible locations of virtual targets from the feasible region, refer to this discrete set as $\mathcal{V}$. The assignment problem is constrained to choose a maximum of $M_V$ virtual targets from $\mathcal{V}$.

The virtual target selection problem reduced to a variant of three dimensional assignment problem on a tripartite graph. Given the set of pursuers, $\mathcal{P}$, set of possible virtual target locations, $\mathcal{V}$, and the set of evaders, $\mathcal{E}$, each pursuer in $\mathcal{P}$ needs to be assigned to a virtual target and an evader. However, unlike the three dimensional assignment problem, not every virtual target needs an assignment. Moreover, the number of virtual targets selected is limited by the operating scenario, and therefore it has to be constrained to a maximum number of virtual targets, $M_V$. Due to these additional constraints, the problem differs from the traditional assignment problem, and needs a new formulation. The cost of assignment is a function of the pursuer ($i$) -virtual target ($k$) - evader ($j$) assignment, as described in the previous section, and denoted as $c_{ijk}$.

The following mixed integer linear program (MILP) formulated this above assignment problem using the binary variables $x_{ijk}$, $i \in \mathcal{P}, j \in \mathcal{E}, k \in \mathcal{V}$:

$$\min \sum_{i \in \mathcal{P}, j \in \mathcal{E}, k \in \mathcal{V}} c_{ijk} x_{ijk}, \tag{25}$$

subject to

$$\sum_{i \in \mathcal{P}, k \in \mathcal{V}} x_{ijk} \geq 1, \forall j \in \mathcal{E}, \tag{26}$$

$$\sum_{j \in \mathcal{E}, k \in \mathcal{V}} x_{ijk} = 1, \forall i \in \mathcal{P}, \tag{27}$$

$$\sum_{i \in \mathcal{P}, j \in \mathcal{E}} x_{ijk} \leq y_k, \forall k \in \mathcal{V}, \tag{28}$$

$$\sum_{k \in \mathcal{V}} y_k \leq M_V, \tag{29}$$

$$x_{ijk}, y_k \in \{0, 1\}, \forall i \in \mathcal{P}, j \in \mathcal{E}, k \in \mathcal{V}. \tag{30}$$

In the above formulation, the cost function in Eq. (25) minimizes the total cost of assignment. The inequality in Eq. (26) enforce the constraint that every evader is assigned to at least one pursuer at a virtual target. The constraint Eq. (27) enforces that every pursuer is assigned to some virtual target and an evader combination. To limit the use of number of virtual targets, we introduce the indicator binary variables $y_k$, $k \in \mathcal{V}$, where a virtual target, $k$, is allowed to be used only if the corresponding $y_k$ is equal to one; this is enforced by inequality Eq. (28). The maximum number of virtual targets that could be used is constrained by equation Eq. (29), and Eq. (30) are the binary constraints.

## V. Example

Consider an example of four pursuers teamed against two evaders taking constant bearing and course. Further, consider that the pursuers must choose three virtual targets so as to minimize the energy for reaching the evaders. In order to demonstrate the approach proposed in this paper, let the initial conditions for the pursuers, evaders, and the domain of the virtual targets be defined in Table 1.

Using the domain for the x and y positions as described in Table 1, the number of virtual target candidates is varied. The combinatorics problem is this: provide a number of candidate virtual targets, locate the optimal locations of a specified number amongst those candidates that minimize the energy for the team of pursuers to reach the team of evaders.

The procedure for conducting this optimization is as follows:

1. For each pursuer-virtual target-evader pair, compute the location and size of each Apollonius circle.

**Table 1 Example Initial Conditions**

| Parameter | X-Position (DU) | Y-Position (DU) |
|---|---|---|
| $P_1$ | 0 | 1 |
| $P_2$ | −1 | 2 |
| $P_3$ | −2 | 6 |
| $P_4$ | −1 | 8 |
| $E_1$ | 20 | 2 |
| $E_2$ | 20 | 8 |
| $VT$ | [3,8] | [−4,10] |

2. For each Apollonius circle, compute the location of the interception point for a pursuer to reach the evader.
3. The cost for each pursuer-virtual target-evader assignment is computed as a function of path lengths and the angle at the virtual target.
4. A linear program is used to evaluate and obtain the optimal assignment of pursuer-virtual target-evader pairs that minimize the overall energy.

In order to conduct the search and optimization, the Julia programming language [42] was used as well as the JuMP package [43] and Cbc (Coin-or branch and cut) open-source mixed integer linear programming solver [44].

Provided the conditions as described in Table 1, the test runs were conducted raging from 9 to 2500 candidate virtual target locations. Some selected results are described in Table 2 and the relationship between run-time (for the specific initial conditions in this problem) are shown in Fig. 2 and the relationship to the optimality for is shown in Fig. 3.

As can be seen in Fig. 2, the convergence time increases with the number of candidate virtual targets.

As can be observed in Fig. 3, the cost for the team of pursuers converges as the number of candidate virtual targets increases. Because the intervals for the virtual targets are held constant; the oscillations in the cost can are directly correlated to candidate value functions laying on the straight-line paths necessary for the pursuers to reach their chosen evader. For example, going from 400 to 2500 candidate locations makes the meshing of the virtual targets much tighter. But in the case where there are 1369 virtual targets, the cost is lower. This is caused by the meshing of the candidate virtual targets aligning with the straight-line paths better than the case where there are 2500 candidates. This leads the authors to believe that a heuristic approach for locating the virtual targets on candidate intersections between pursuer-evader paths may lead to faster and more accurate virtual target locations when there are fewer evaders than pursuers. While Fig. 3 allows the reader to draw this conclusion, the heuristic approach is left as future work.

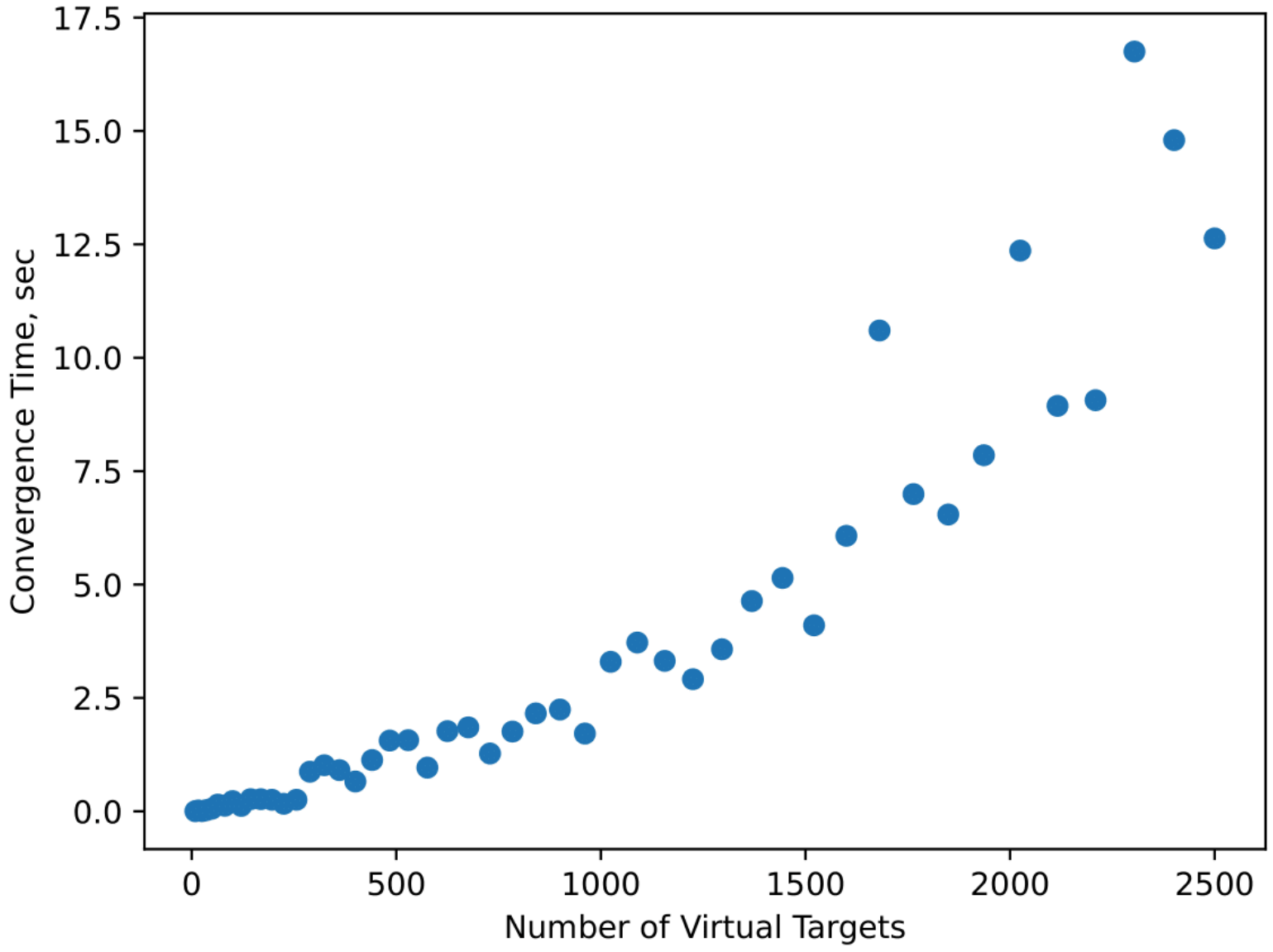

**Fig. 2 The relation between the convergence time with respect to the number of candidate virtual targets**

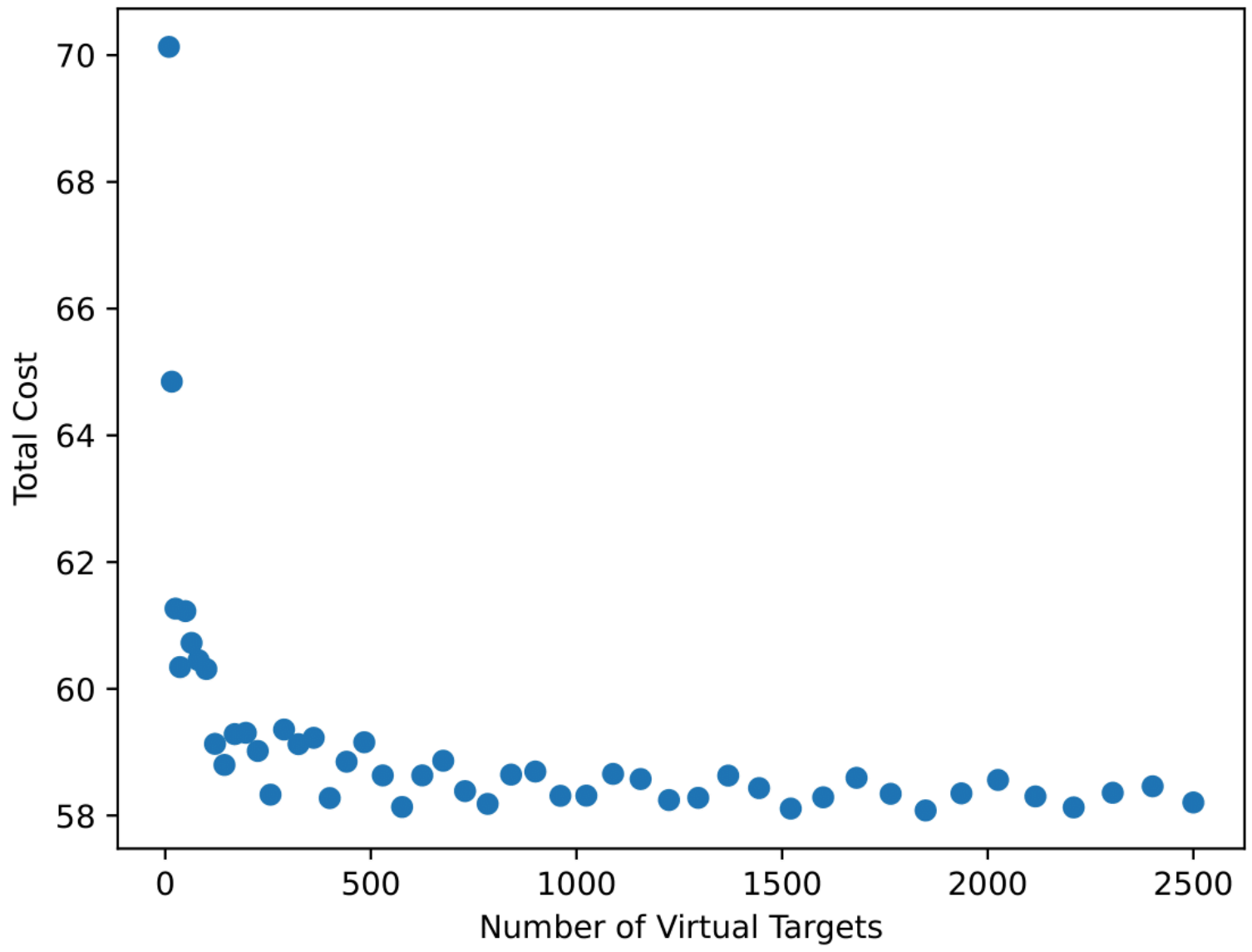

**Fig. 3 The value of the optimization converges with regard to the number of candidate virtual targets**

Consider the case where there are 100 candidate virtual targets, all shown in Fig. 4. In this case, one can observe the blue grid of points that make up the candidate virtual target locations in the left figure. In the right figure, the optimal virtual targets are selected that minimize the objective in Eq. (1). This results in almost straight-line behavior for each of the pursuers to reach their targets. In the left figure, every single Apollonius circle is shown, this represents the locus of possible interceptions of pursuers to each evader provided they stop at each possible virtual target along the way.

Taking a closer investigation as the number of virtual targets is seen in Fig. 5.

From Fig. 5, one can observe that when a very small number of targets is selected, a large maneuver is demonstrated for nearly every pursuer to reach the evaders, only the bottom-most pursuer reaches the evader without needing to maneuver at one of the candidate virtual targets. However, as the number of candidate virtual targets increases to 16 or even 400, the paths begin to straighten out, thereby reducing the objective cost because the maneuver at the virtual target is minimized. Ultimately, when a very fine grid is selected, say 2500 potential virtual target locations, the paths are straight-line paths. While this observation is true for this example; it can be stated that as the number of candidate virtual points increases, the paths will get straighter and shorter; this is reflected in the cost being reduced (as highlighted in Fig. 3).

**Table 2 Selected Experimental Runs**

| Number of Candidate $VT$s | Time (sec) | Objective Value |
|---|---|---|
| 9 | 0.002556 | 70.132 |
| 16 | 0.00245 | 64.849 |
| 25 | 0.00555 | 51.265 |
| 100 | 0.2510 | 61.311 |
| 400 | 0.6712 | 58.275 |
| 2500 | 11.252 | 58.205 |

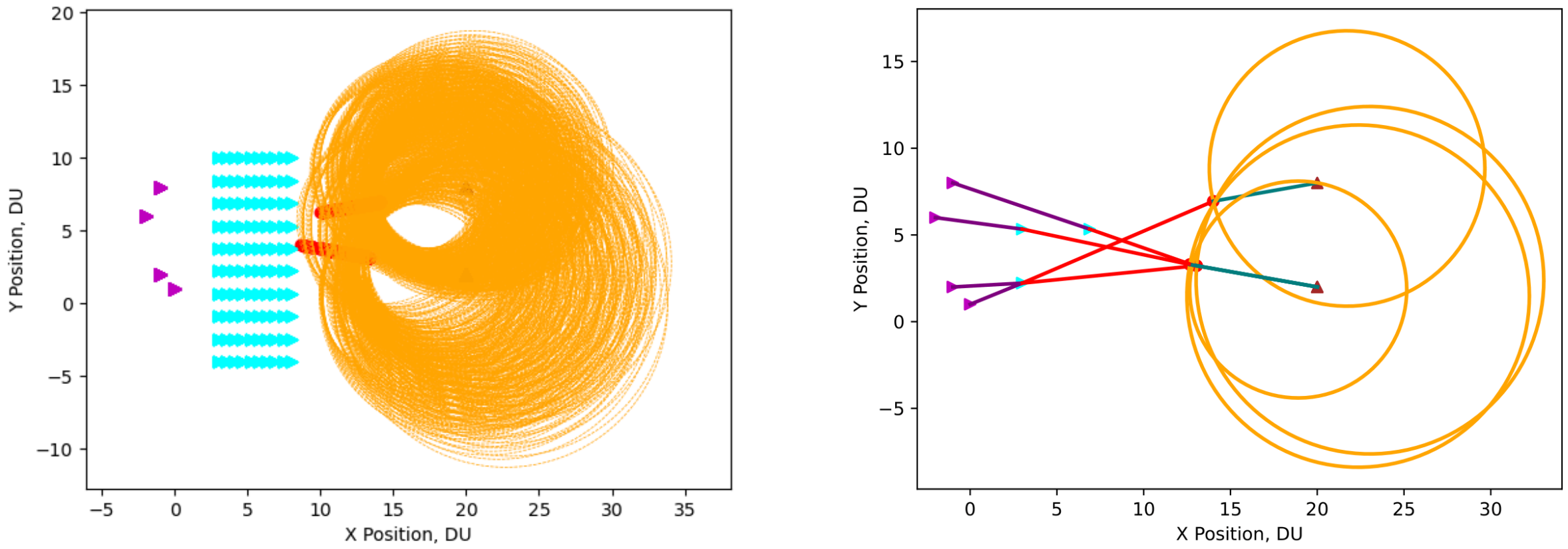

**Fig. 4 100 Sample Virtual Target Locations**

Straight-line paths may not always be feasible by the pursuers to reach each of the evaders. This could happen when the number of pursuers and evaders greatly out-numbers the number of potential virtual target locations. Finding the critical conditions for which straight line paths are not feasible is outside the scope of this work and is left for future investigation.

## VI. Conclusions

In conclusion, a multi-pursuer multi-evader assignment problem has been posed wherein the pursuers must first navigate to a set of virtual targets, where they assign roles, and then pursue a team of evaders. Using a Apollonius geometry and a linear program solver; the virtual targets are obtained that minimize the path length of each of the evaders as well as the maneuver made at each of the virtual targets – the energy of the pursuing team is minimized.

Future work could consider a number of interesting research directions including, but not limited to, 1) Ensuring path deconfliction of solutions, 2) simulating higher-fidelity vehicles with more realistic performance constraints and measurement uncertainty 3) Performing software-in-the-loop and hardware-in-the-loop tests.

In order to ensure path deconfliction between agents; one could impose a kinematic constraint as an additional constraint in the mixed integer linear program. Such a constraint would take a given assignment and then require that for that assignment that the miss distance between agents be lower bounded. This non-trivial change may provide utility when ranges are small between pursuers, evaders, and virtual targets. To consider higher fidelity models or the introduction of measurement models can inform robust strategies for weapon-target-assignment. These approaches would require more computational effort compared to the approach provided in this paper; but, would add credence to the solutions and viability for future hardware testing. Lastly, any amount of software or hardware integration would serve useful, as it would highlight any discrepancies or considerations that weren't accounted for in the current formulation. A proper simulation may also provide the same utility; but, the presence of external disturbances, communication delays, measurement noise, and a myriad of real-world phenomena could be easier to test than to simulate outright.

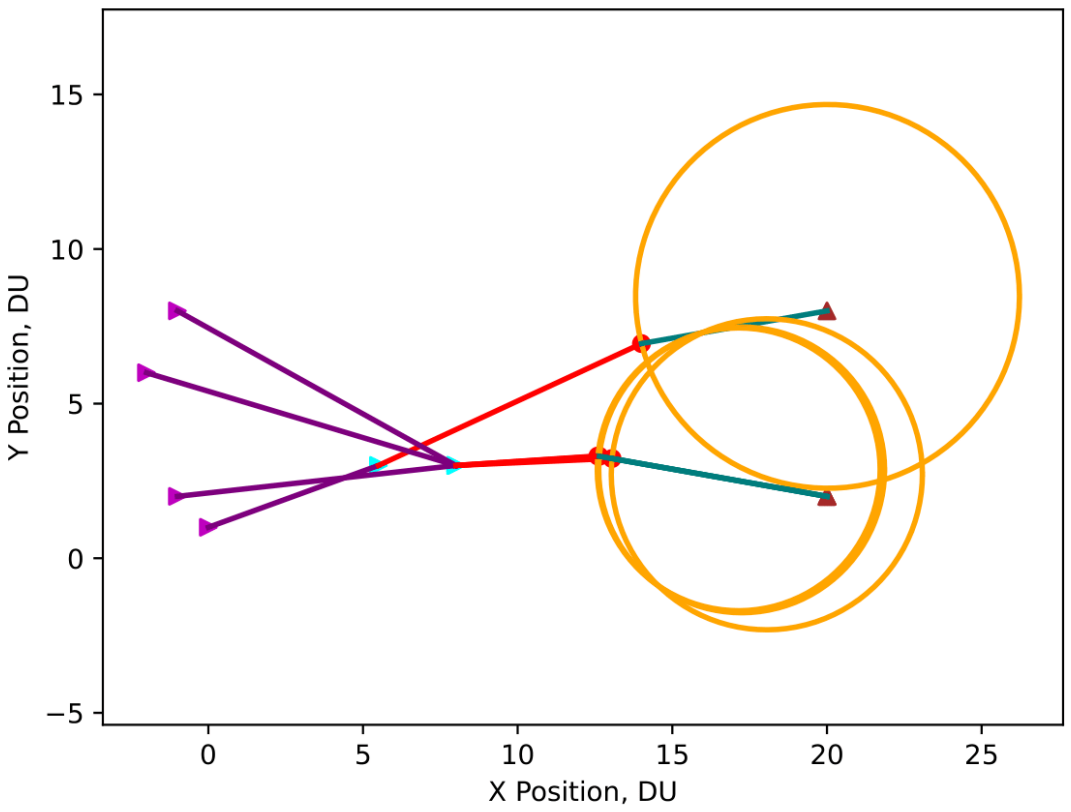

3/9 Virtual Targets

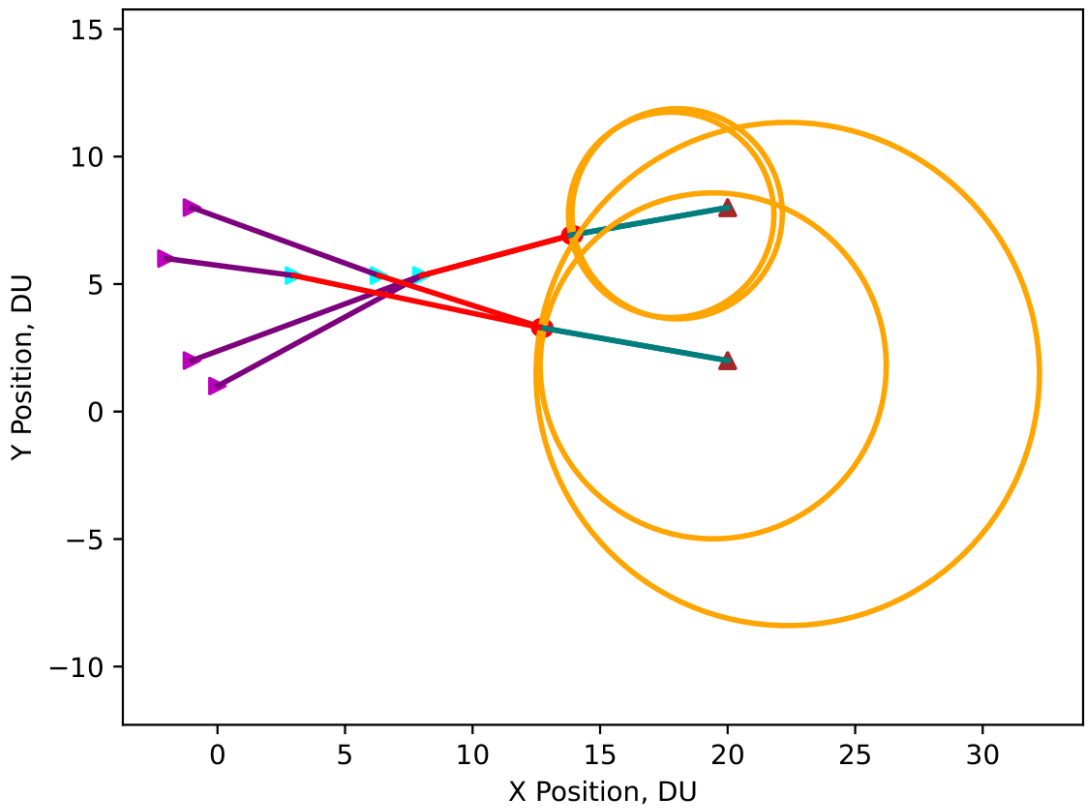

3/16 Virtual Targets

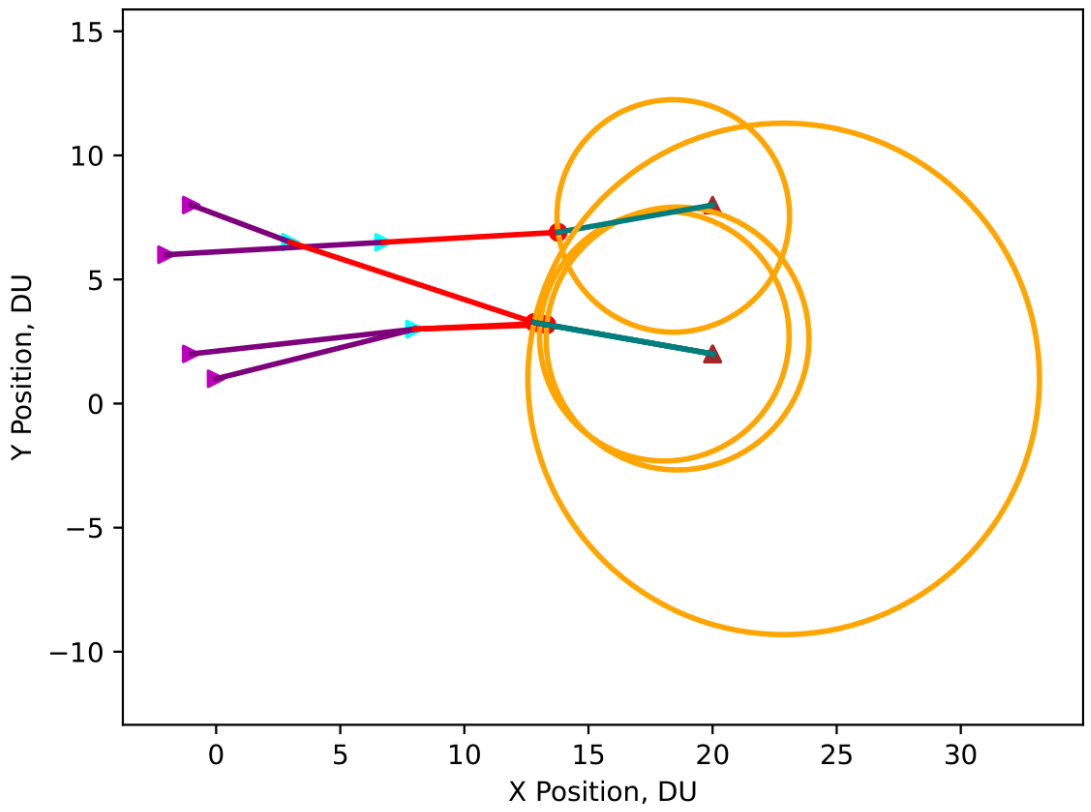

3/25 Virtual Targets

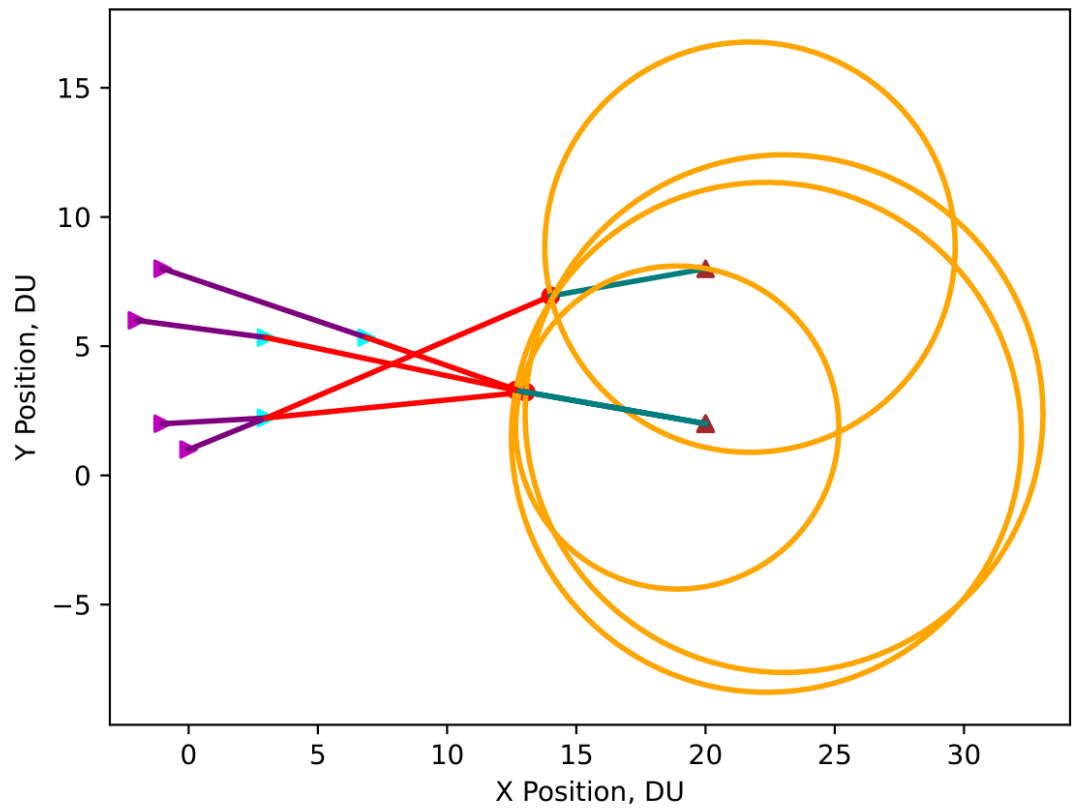

3/100 Virtual Targets

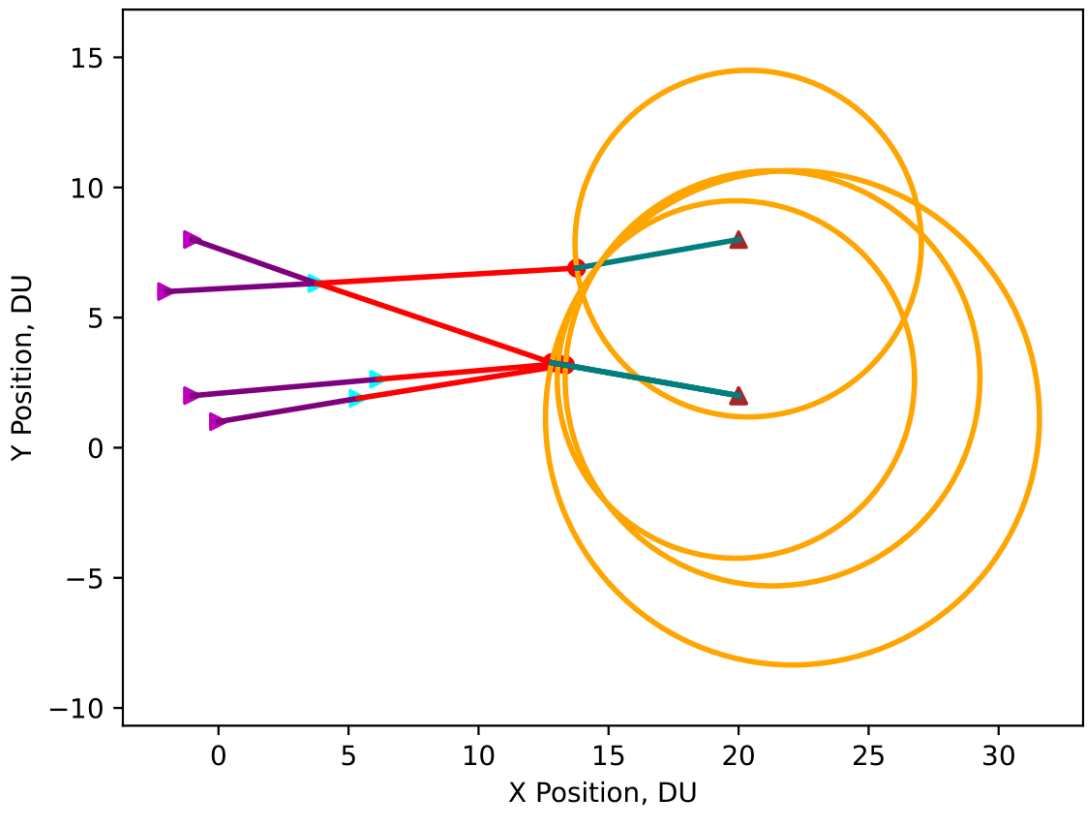

3/400 Virtual Targets

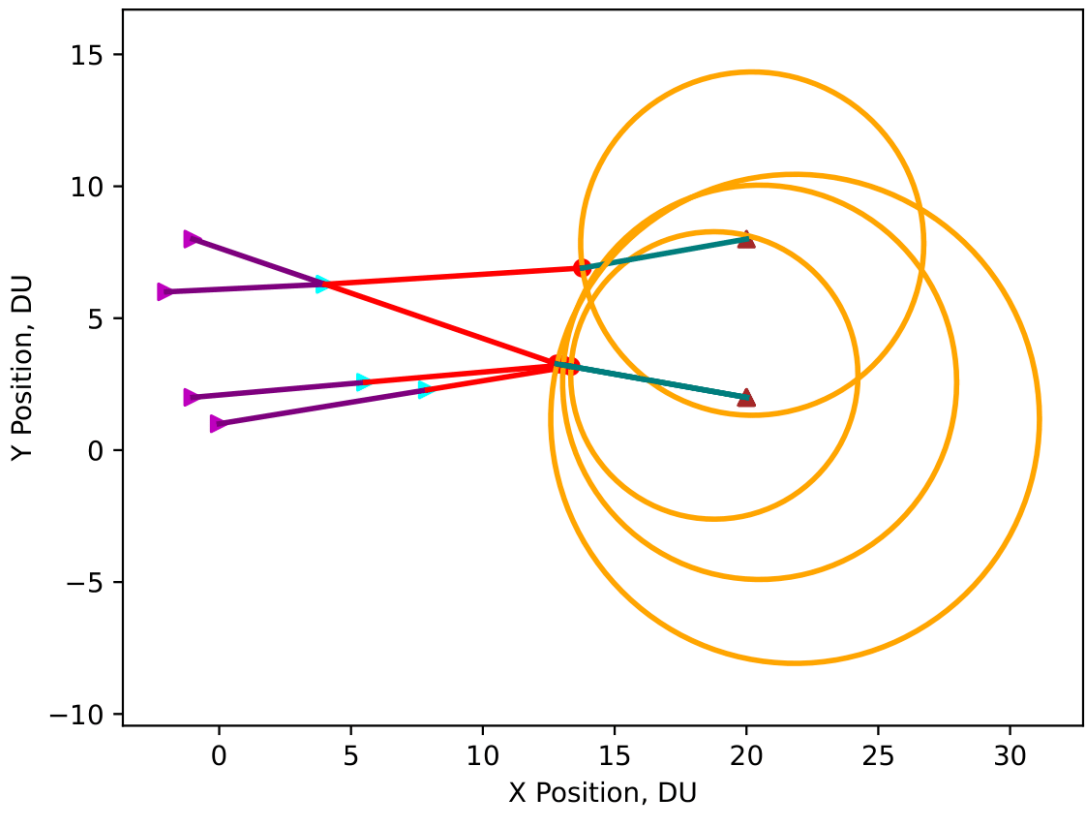

3/2500 Virtual Targets

**Fig. 5 Optimized solutions for the location of virtual targets where $\frac{a}{b}$ virtual targets represent $a$ selected number of virtual targets for $b$ candidates over the field $VT_x \in [3, 8]$ and $VT_y = \in [-4, 10]$.**

## VII. Acknowledgments

The authors would like to thank Dr. Adam Gerlach for providing assistance with using the Julia programming language. DISTRIBUTION STATEMENT A. Approved for public release: distribution unlimited (AFRL-2023-2613).